\documentclass[reqno]{amsart}
\usepackage{amssymb}

\newcommand{\rankg}{r}
\newcommand{\ranky}{{n}}
\newcommand{\rankh}{{\ell}}
\newcommand{\C}{{{\mathbb C}}}
\newcommand{\nnumber}{{\mathbb N}}
\newcommand{\integer}{{\mathbb Z}}
\newcommand{\nonnegative}{{{\mathbb Z}_{\geq 0}}}
\newcommand{\Sym}[1]{{\mathfrak{S}}_{#1}}
\newcommand{\dual}[1]{#1^{\ast}}
\newcommand{\tensorN}{^{\tensor \rankh}}
\newcommand{\uesita}[2]{{\substack{#1 \\ #2}}}
\newcommand{\Center}{{ Z}}
\newcommand{\lie}[1]{{\mathfrak {gl}}_{#1}}

\newcommand{\Mat}[2]{E_{#1}}
\newcommand{\BGG}{{\mathcal O}_{\rankg}}
\newcommand{\nil}[2]{{\mathfrak n}_{#2}^{#1}}

\newcommand{\Id}[1]{{ I}_{#1}}
\newcommand{\tensorl}{^{\otimes \rankh}}

\newcommand{\dominant}[1]{{P}_{{#1}}^{+}}
\newcommand{\indominant}[1]{P_{#1,\integer}^+}
\newcommand{\Hcat}{{\mathcal C}_{\Hecke{\rankh}}}
\newcommand{\YNcat}{{\mathcal C}_{\YN}}
\newcommand{\YSNcat}{{\mathcal C}_{\YSN}}

\DeclareMathOperator{\End}{End}

\DeclareMathOperator{\Wt}{Wt}

\DeclareMathOperator{\sign}{sign}
\DeclareMathOperator{\gr}{gr}

\newcommand{\YN}{{ Y}({\mathfrak{gl}}_\ranky)}
\newcommand{\YSN}{{ Y}({\mathfrak {sl}}_\ranky)}
\newcommand{\Vbox}{{\C^\ranky}}
\newcommand{\tensor}{\otimes}
\newcommand{\UN}{{ U}({\mathfrak {gl}}_\ranky)}
\newcommand{\qdet}{{\rm qdet}\,}

\newcommand{\ev}[1]{{{\rm {ev}}_{#1}}}
\newcommand{\Ystandard}{{ M}}

\newcommand{\Ysimple}{{ V}}

\newcommand{\Dfun}[1]{D_{#1}}

\newcommand{\functorH}{{ F}_{\lambda}}

\newcommand{\cartan}[1]{{\mathfrak h}_{#1}}
\newcommand{\e}[1]{{\epsilon_{#1}}}
\newcommand{\Root}[1]{R_{#1}}
\newcommand{\Hecke}[1]{{{\mathcal H}_{#1} }}

\newcommand{\standard}{{\mathcal K}}
\newcommand{\simple}{{\mathcal L}}
\newcommand{\Hmodule}{{ M}}

\title{
Drinfeld Functor
and Finite-Dimensional Representations
of Yangian }
\author{Tomoyuki Arakawa }
\address{Graduate school of Mathematics,  Nagoya University,
Chikusa-ku, Nagoya, 464-8602, JAPAN}
\email{tarakawa@math.nagoya-u.ac.jp}
\begin{document}
\theoremstyle{plain}
  \newtheorem{Th}{Theorem}
  \newtheorem{Pro}[Th]{Proposition}
  \newtheorem{Lem}[Th]{Lemma}
   \newtheorem{Facts}[Th]{Facts}
\theoremstyle{definition}
  \newtheorem{Def}[Th]{Definition}
  \newtheorem{Co}[Th]{Corollary}
  \newtheorem{dfandpr}[Th]{Definition and Proposition}
\theoremstyle{remark}
 \newtheorem{Rem}[Th]{Remark}
 \newtheorem{Ex}[Th]{Example}
\bibliographystyle{alpha}
\numberwithin{equation}{subsection}

\begin{abstract}
{We extend the results of 
Drinfeld on the Drinfeld functor
to the case $\rankh\geq\ranky$.
We present the character of 
finite-dimensional representations
of the Yangian $\YSN$ in terms
of the Kazhdan-Lusztig polynomials
as a consequence.}
\end{abstract}

\maketitle

\section*{Introduction}
In this article we study the representations
of  the Yangian $\YSN$.
The Yangian is a quantum group introduced
 by V. G. Drinfeld (\cite{D3}).
The parameterization of the simple finite-dimensional representations
of $\YSN$ was obtained in \cite{D1}
by the sequences of monic polynomials $Q(u)=
(Q_1(u),\dots,Q_{\ranky-1}(u))$
called the Drinfeld polynomials.
Furthermore,
he has constructed 
in 
 \cite{D2} a  functor 
from the category 
$\Hcat$ of finite-dimensional representations
of
 the degenerate affine Hecke algebra $\Hecke{\rankh}$
to
the category
$\YSNcat$ of finite-dimensional representations
of
$\YSN$.
This functor is called  the Drinfeld functor.
It was stated in \cite{D2} that
 as well as the classical Frobenius-Schur duality,
 the Drinfeld functor gives
the categorical equivalence between $\Hcat$ and the certain
subcategory of $\YSNcat$ in the case
$\rankh<\ranky$.
 Chari-Pressley 
generalized this
 duality
to the case
between the affine Hecke algebra
and the quantum affine algebra.
They proved that the categorical equivalence holds
in this case
 as well
provided that $\rankh<\ranky$
(\cite{CP1}).

However, due to the restriction $\rankh<\ranky$, 
the above categorical equivalence does not describe all the finite-dimensional
representations of the Yangian $\YSN$.
In particular, even the characters of finite-dimensional
representations of
$\YSN$ have  not been known, except for
the case  $\ranky=2$ (\cite{CP3}) and 
the special class of the representations called tame
(\cite{Tara2}).

The main purpose of this article is to extend
the Drinfeld's results \cite{D2} to the case $\rankh\geq \ranky$.
To be more precise,
 we  first show  the followings 
without restriction $\rankh<\ranky$:
\begin{enumerate}
\item
\label{results:standard}
 The Drinfeld functor sends  the standard modules of $\Hecke{\rankh}$
to  zero or the  highest modules of $\YSN$
(Theorem \ref{Th:properties_of_Dfun_standard}).
\item  
\label{results:simple}
The Drinfeld functor sends 
the simple modules of $\Hecke{\rankh}$
to zero or the simple modules of $\YSN$
(Theorem \ref{Th:properties_of_Dfun_simple}).
\end{enumerate}
Here the standard modules 
are certain induced $\Hecke{\rankh}$-modules
which have unique simple quotients
(see subsection \ref{subsection:Hecke_rep}).
We also determine the 
explicit images of the standard modules.
It turns out that the highest weight modules obtained
as the images of the standard modules are those tensor product
modules of the evaluation representations
studied in  \cite{Ak}.
We note that any simple $\YSN$-module
is isomorphic to the image of a simple
$\Hecke{\rankh}$-module for some $\rankh$.

Further, combining the above results
with that of
the representation theory 
of $\Hecke{\rankh}$,
we state the following:
\begin{enumerate}
\addtocounter{enumi}{2}
\item The multiplicity
formula of $\YSN$  expressed in terms of
the Kazhdan-Lusztig polynomials  (Theorem \ref{th:main})\footnote{
After completing this article,
the author was notified that
E. Vasserot obtained the similar formula
in terms of 
intersection cohomologies
in the case of the quantum affine algebra
by geometrical method
([preprint, math.QA/9803024]). }.
\end{enumerate}
This is the result of
considering the composition
$\Dfun{\rankh}\circ F_{\lambda}$
of the two exact functors,
where  $F_{\lambda}$ is the functor 
 from the Bernstein-Gelfand-Gelfand
category $\BGG$ of 
the complex Lie algebra
$\lie{\rankg}$ to the category $\Hcat$,
obtained by Suzuki and the author in \cite{AS}.

\medskip



\bigskip

\section{Preliminaries}
\subsection{Yangian}
Let $\ranky$ be a positive integer.
First we review some
fundamental facts about the algebra structure of
the {\emph {Yangian}} $\YSN$.
Our main references are
\cite{D3,D1,MNO} and 
we basically follow the notation
of \cite{MNO}.

\smallskip
Let 
 $$R(u)=1-\dfrac{P}{u}\in \End(\C^\ranky\otimes \C^\ranky),$$
where $P$ is the permutation operator in 
$\Vbox\tensor \Vbox$ and $u$ is a parameter.
Let $\Mat{ij}{\ranky}\in \lie{\ranky}=\lie{\ranky}(\C)$
denote the usual matrix operator on $\C^{\ranky}$.
The Yangian $\YN$ is the unital associative algebra over $\C$
with generators $t_{ij}^{(k)}$ ($1\leq i,j\leq \ranky$, 
$k=0,1,2,\dots $) and the defining relations
\begin{eqnarray}
\label{ternary}
R_{12}(u-v)t_1(u)t_2(v)=t_2(v) t_1(u)R_{12}(u-v),
\end{eqnarray}
where 
\begin{equation*}
\begin{array}{l}
t(u)=\sum_{i,j} t_{ij}(u)\tensor \Mat{ij}{\ranky}\in \YN
\tensor \End(\Vbox),\\
~~t_{ij}(u)=\sum_{d=0}^{\infty}t_{ij}^{(d-1)}u^{-d}
\in \YN [[u^{-1}]]
\end{array}
\end{equation*}
Here we put 
$t_{ij}^{(-1)}=\delta_{ij}{\rm id}$ and
both sides of  (\ref{ternary}) are regarded as
elements of $\YN((v^{-1}))[[u^{-1}]]\otimes \End(\Vbox)\otimes \End(\Vbox)$ and
the subindexes of  $t(u)$ 
and $R(u)$
indicate to which copy of $\End(\Vbox)$
these matrices correspond.

The defining relations \eqref{ternary} are equivalent to the following 
relations:
\begin{align}
\label{eq:y_def2}
[ t_{ij}^{(r)},t_{kl}^{(s-1)}]
-[ t_{ij}^{(r-1)},t_{kl}^{(s)}]
=t_{kj}^{(r-1)}t_{il}^{(s-1)}-
t_{kj}^{(s-1)}t_{il}^{(r-1)}
\end{align}
($1\leq i,\,j,\,k,\,l\leq \ranky$,
$r,s\in\nonnegative$) (\cite{MNO}).

The algebra
$\YN$ is a Hopf algebra with coproduct
\begin{equation}
\label{eq:coproduct}
\Delta:\, t_{ij}(u)\mapsto
\sum_{a=1}^{\ranky} t_{ia}(u)\tensor t_{aj}(u),
\end{equation}
antipode $S(t(u))=t(u)^{-1}$ and counit 
$\varepsilon (t(u))=1$.

Let $U(\lie{\ranky})$
denote the universal enveloping algebra 
of 
the  Lie algebra 
$\lie{\ranky}$.
The algebra $U(\lie{\ranky})$
is considered as a subalgebra
of $\YN$ 
by the inclusion homomorphism
defined by
\begin{equation*}
\begin{array}{ccc}
\UN &\longrightarrow & \YN\\
\Mat{ij}{\ranky}&\longmapsto & t_{ij}^{(0)}.
\end{array}
\end{equation*}
On the other hand, for $a\in \C$, the map
\begin{equation}
\begin{array}{lccc}
\label{evaluation}
{\rm ev}_a:&\YN &\longrightarrow & \UN\\
&t_{ij}(u)&\longmapsto &
\delta_{ij}+\dfrac{\Mat{ij}{\ranky}}{ u-a}
\end{array}
\end{equation}
defines an algebra homomorphism.
For a $\lie{\ranky}$-module $V$,
let $\ev{a}(V)$ denote the $\YN$-module
 obtained by pulling $V$ back by \eqref{evaluation}.

The {\emph {quantum determinant}} $\qdet t(u)$
is defined as
$$ \qdet t(u) =\sum_{w\in W_{\ranky}} (-1)^{{\rm sgn}(w)}
t_{w(1),1}(u)t_{w(2),2}(u-1)\cdots t_{w(\ranky),\ranky}(u-\ranky+1),$$
where $W_{\ranky}$ denotes the 
symmetric group $\Sym{\ranky}$.
The coefficients of the quantum determinant are algebraically
independent and generate the center $\Center(\YN)$
of $\YN$.

For a formal series $f(u)=1+f_1 u^{-1}+f_2 u^{-2}+\dots
\in \C[[u^{-1}]]$, the multiplication
\begin{equation}
\label{auto1}
t(u)\longmapsto f(u)t(u)
\end{equation}
defines an automorphism of $\YN$.
It is known that the Yangian $\YSN$ can be defined as the subalgebra of $\YN$
consisting of elements fixed by all automorphisms of the form
(\ref{auto1}) (\cite{MNO}). One has 
a tensor product decomposition
\begin{equation}
\label{Y-center}
\YN\cong { Z}(\YN)\tensor \YSN.
\end{equation}
Hence any $\YN$-module can be considered as a $\YSN$-module.

\medskip
Let
\begin{align}
\label{eq:filterYN}
F_i \YN:=
\sum_\uesita{k\in\nonnegative}{1\leq i_a,j_a\leq \ranky}
\sum_\uesita{r_i\in \nonnegative}{r_1+\dots+r_k\leq i}
 \C t_{i_1 j_1}^{(r_1)}\dots
t_{i_k j_k}^{(r_k)}\quad
\left( i\in\nonnegative\right).
\end{align}
Then, by \eqref{eq:y_def2},
$F_i \YN\cdot F_j \YN\subset
F_{i+j}\YN$,
and
\eqref{eq:filterYN}
defines a filteration on $\YN$.
Let $\gr \YN$ denote
 the corresponding the
graded algebra;
 $$\gr \YN=\bigoplus_{i\in\nonnegative}F_i \YN
/ F_{i-1}\YN\quad
\left( F_{-1} \YN=0\right).$$
Then,
$\gr \YN$ is isomorphic to $U(\lie{\ranky}[t])$,
where  $U(\lie{\ranky}[t])$
is the universal enveloping algebra of
the polynomial current algebra
$\lie{\ranky}[t]:=\lie{\ranky}\tensor \C[t]$  
with $\nonnegative$-grading such that
the degree of the element $X\tensor t^r$
$\left( X\in \lie{\ranky}\right)$
equals $r$ (\cite{MNO}).

\subsection{Drinfeld polynomials}
In this subsection we review the classification
theory of finite-dimensional simple $\YSN$-modules studied by Drinfeld
(\cite{D1}, see also \cite{CP2,Mo})

A representation $V$ of $\YN$ is called 
{\emph {highest weight}}
if there exists a cyclic vector $v$ such that
$t_{ij}(u) \cdot v=0$ ($1\leq i<j\leq \ranky$) and 
$t_{ii}(u-i)\cdot v= \zeta_i(u) v$ ($1\leq i\leq \ranky$)
for some formal series $\zeta_{i}(u)\in 
1+u^{-1}\C[[u^{-1}]]$. 
The vector $v$ is called the {\emph {highest weight vector}}
of $V$
and the sequence $\zeta(u)=(\zeta_1(u),\dots,\zeta_\ranky(u))$
is called the
{\emph {highest weight }}of $V$.
The central element $\qdet t(u)$   acts as
a constant
$\zeta_1(u)\zeta_2(u-1)\dots,\zeta_\ranky(u-\ranky+1)$
on a highest weight module $V$.
As in the classical Lie algebra theory,
any highest weight $\YN$-module has a unique simple quotient,
in which the image of its highest weight vector is nonzero. 

It is known by Drinfeld that a simple highest weight module of $\YN$ 
is finite-dimensional if and only if there exists a sequence
of monic polynomials
$Q(u)=\left(Q_1(u),\dots,Q_{\ranky-1}(u)\right)$ such that 
\begin{equation}
\label{df:D_poly}
\frac{\zeta_{k}(u) }{\zeta_{k+1}(u+1)}=
\frac{Q_k(u+1)}{ Q_{k}(u)}
\end{equation}
for $k=1,\dots,\ranky-1$.
A theorem of Drinfeld states  that there is a one-to-one correspondence
between the finite-dimensional simple
 $\YSN$-modules
and   the sequences of monic polynomials
$Q(u)=(Q_1(u),\dots,Q_{\ranky-1}(u))$
defined by  \eqref{df:D_poly} (\cite{D1}).
 The $Q(u)$ are called
{\emph {Drinfeld polynomials}}.
\begin{Rem}
The standard symbol for the Drinfeld polynomials
is $P(u)$.
However, 
this symbol for a polynomial is reserved
for the Kazhdan-Lusztig polynomials in this article.
\end{Rem}

\subsection{Degenerate Affine Hecke Algebra}
\label{degenerate Hecke}
Let $\rankh$ be a positive integer.
Let ${\cartan{\rankh}}$
be the Cartan subalgebra of
$\lie{\rankh}$, which consists of the diagonal matrices.
Define a basis $\{\e{i}\}_{i=1}^{\rankh}$
 of $\cartan{\rankh}$ by putting $\e{i}=E_{ii}$.
The dual space $\dual{\cartan{\rankh}}$ is identified
with $\cartan{\rankh}$ via the inner product 
$\left<\e{i}, \e{j}\right>=\delta_{ij}$.

Let $\Root{\rankh}$ be the root system of $\lie{\rankh}$;
$$ 
\begin{array}{l}
\Root{\rankh}=\left\{\alpha_{ij}=\e{i}-\e{j}\mid
1\leq i\ne j\leq \rankh\right\},\\
\Root{{\rankh}}^+=\left\{\alpha_{ij}\in \Root{\rankh}\mid
i<j \right\},
\\
\Pi_{\rankh}=\left\{\alpha_i=\alpha_{ii+1}\mid
i=1,\dots,\rankh-1\right\},
\end{array}
$$
where $\Root{\rankh}^+$ is the set of positive roots
and
$\Pi_{\rankh}$ is the set of simple roots.
Let $
\rho=\frac{1}{2}\sum_{\alpha\in R_{\rankh}^+}\alpha$.
We identify the coroots with the roots throughout 
this article.

Let $s_{\alpha}\in W_{\rankh}$ denote the reflection
corresponding to $\alpha\in R_{\rankh}$;
$$s_{\alpha}\cdot\lambda=\lambda-
\lambda(\alpha)\alpha\quad(\,\lambda\in\dual{\cartan{\rankh}}\,).$$
Put $s_{ij}=s_{\alpha_{ij}}$ ($\alpha_{ij}\in R_{\rankh}$)
and $s_i=s_{\alpha_{i}}$ ($\alpha_i\in \Pi_{\rankh}$).

Let $S(\cartan{{\rankh}})$ be the symmetric algebra of 
$\cartan{{\rankh}}$, which is isomorphic to the
polynomial ring over $\dual{\cartan{{\rankh}}}$.

The 
{\emph {degenerate affine Hecke algebra}}
 $\Hecke{{\rankh}}$ of $GL_{\rankh}$ (\cite{D2})
is the associative algebra over $\C$
such that
$$\Hecke{{\rankh}}\cong \C[W_{\rankh}]\tensor S(\cartan{{\rankh}})$$
as a vector space, the subspaces
$\C[W_{{\rankh}}]\tensor \C$ and $\C\tensor S(\cartan{{\rankh}})$
are subalgebras of $\Hecke{{\rankh}}$,
and the following relations hold in it;
\begin{equation}
\label{eq:h_def}
s_{\alpha}\cdot {\xi}
-{s_{\alpha}(\xi)}\cdot s_{\alpha}
=-\alpha(\xi) ~~~( \,\alpha\in \Pi_{{\rankh}}, \,
\xi\in \cartan{{\rankh}}\,),
\end{equation}
where the elements
${\xi}\in\cartan{{\rankh}}$ and $w\in W_{\rankh}$
are identified with $1\tensor \xi\in\Hecke{{\rankh}}$
and $w\tensor 1\in\Hecke{\rankh}$ respectively.
One has 
\begin{align}
\label{eq:Hecke_dec}
\Hecke{\rankh}=\C[W_{\rankh}]\cdot S(\cartan{{\rankh}})=
S(\cartan{{\rankh}})\cdot \C[W_{\rankh}].
\end{align}
We put $\Hecke{0}=\C$ for convenience.

The center $\Center(\Hecke{{\rankh}})$ of this algebra
equals the $W_{{\rankh}}$-invariant polynomials
$S(\cartan{{\rankh}})^{W_{{\rankh}}}=
\C[\e{1},\dots,\e{{\rankh}}]^{W_{{\rankh}}}$
(\cite{Lus_Hecke}).

Define elements $y_i\in \Hecke{\rankh}$
$\left(i=1,\dots,\rankh\right)$
by
$$y_i:=s_{1i}\cdot \e{i}\cdot s_{1i}=
\e{i}-\sum_{j<i}s_{ji}.$$
Then one can see by direct calculations that
\begin{align}
&w\cdot y_i=y_{w(i)}\cdot w\quad \left(\,i=1,\dots,\rankh,\,w\in
W_{\rankh}\right),\label{eq:another_gen1}\\
&[y_i,y_j]=-(y_i-y_j)s_{ij}
\quad \left(1\leq i,j\leq \rankh\right).
\label{eq:another_gen2}
\end{align}
The algebra $\Hecke{\rankh}$
is isomorphic to the $\C$-algebra
with generators
$w\in W_{\rankh}$
and $y_i$ $\left(i=1,\dots,\rankh\right)$
with
the defining relations
\eqref{eq:another_gen1}, \eqref{eq:another_gen2} and the 
Coxeter relations
of $w$'s
in $W_{\rankh}$.

Let
\begin{align}
\label{eq:filter_H} 
F_i \Hecke{\rankh}:=
\sum_{w\in W_{\rankh}}
\sum_{d_1+d_2+\dots+d_{\rankh}\leq i}
\C  y_1^{d_1}y_2^{d_2}\dots y_{\rankh}^{d_{\rankh}}
w
\quad \left(
i\in\nonnegative\right).
\end{align}
Then,
by \eqref{eq:another_gen1}, 
$F_i \Hecke{\rankh}\cdot F_{j}\Hecke{\rankh}\subset
F_{i+j}\Hecke{\rankh}$,
and
\eqref{eq:filter_H} defines a filteration
on 
 $\Hecke{\rankh}$.
Let $\gr \Hecke{\rankh}$
denote the corresponding graded algebra;
$$\gr \Hecke{\rankh}=
\bigoplus_{i\in \nonnegative} F_i \Hecke{\rankh}/
F_{i-1} \Hecke{\rankh}\quad
\left( F_{-1}\Hecke{\rankh}=0\right).$$
Let $\bar y_i$ denote the image of $y_i$ in
$\gr \Hecke{\rankh}$.
Then 
$\gr \Hecke{\rankh}$
is isomorphic to the graded $\C$-algebra
generated by
$\C[W_{\rankh}]$
and
the polynomial ring
$
\C[\bar y_1,\dots,\bar y_{\rankh}]$
with the relations
 $w\cdot \bar y_i=\bar y_{w(i)}\cdot w$
$\left(i=1,\dots,\rankh,\,w\in W_{\rankh}\right)$,
whose 
 grading is given by $\deg(y_i)=1$ and $\deg(w)=0$.
In particular,
$\gr \Hecke{\rankh}\cong
\C[\bar y_1,\dots,\bar y_{\rankh}]
\tensor_{\C}\C[W_{\rankh}]
$ as a $\C$-vector space.


\medskip
\subsection{Representations of degenerate affine Hecke algebra}
\label{subsection:Hecke_rep}
In this subsection
we review the theory of
the representations of 
$\Hecke{\rankh}$ 
studied in \cite{Ze1,Ro,Gi} (see also \cite{Gi2}),
along the line
introduced in \cite{AS} and developed in \cite{S}.

Let $\Hcat$ denote
the category of  finite-dimensional representations
of $\Hecke{\rankh}$.
Let $\rankg$ be a nonnegative integer.
The representation theory
of $\Hecke{\rankh}$ is well-described
by the language of the Bernstein-Gelfand-Gelfand category $\BGG$
of $\lie{\rankg}$,
via the functors $F_{\lambda}:
\BGG\rightarrow \Hcat$ 
constructed in \cite{AS,S}:

Let
\begin{align*}
&\dominant{\rankg}:=
\left\{
\lambda\in \dual{\cartan{\rankg}}\mid
\lambda(\alpha)\not\in -1,-2,\dots\text{ for all }
\alpha\in R_{\rankg}^+\right\},\\
&\indominant {\rankg}:=\bigoplus_{i=1}^{\rankg}\integer \e{i}
\subset \dominant{\rankg}.
\end{align*}
An element of $\dominant{\rankg}$ (resp. $\indominant{\rankg}$) is called the
dominant
(resp. integral dominant) weight.

For $\lambda\in\dual{\cartan{\rankg}}$,
there is a functor
form $\BGG$
to $\Hcat$
defined by
\begin{align}
\label{functor:H}
\functorH(X)
:&=H_0(\nil{-}{\rankg},X\otimes (\C^{\rankg})\tensorl)_{\lambda-\rho}
\quad (\,X\in \BGG\,)
\\&
=[X\otimes (\C^{\rankg})\tensorl/
\nil{-}{\rankg}(X\otimes (\C^{\rankg})\tensorl)]_{\lambda 
-\rho},
\nonumber
\end{align}
where
$\nil{-}{\rankg}$ denotes the
nilpotent 
subalgebra of $\lie{\rankg}$
generated by the lower triangular
matrices $
\left<\Mat{ij}{\rankg}|i-j>0 \right>$
and
 $X_{\lambda}$ denotes the weight
space of weight $\lambda$ of a $\lie{\rankg}$-module $X$.
The  action of $\Hecke{{\rankh}}$
on the space $\functorH(X)$
is given by
\begin{align*}
&\e{i}\longmapsto \sum_{j=0}^{i} \Omega_{ji}+\frac{\rankg-1}{2}
\quad\left(1\leq i\leq {\rankh}\right)\\
&s_i\longmapsto \Omega_{i i+1}
\quad\left(1\leq i\leq {\rankh}-1\right),
\end{align*}
where
$\Omega_{ij}$
denotes
an endomorphism of $X\otimes (\C^{\rankg})\tensorl$
which acts as 
the Casimir $\Omega
=
\sum_{rs}\Mat{rs}{\rankg}\tensor
\Mat{sr}{\rankg}$
on the tensor product of $i$-th and $j$-th
factors and by identity on all the  other factors.
Here the $0$-th factor corresponds to $X\in\BGG$.
The functor $F_{\lambda}$ is exact if
$\lambda\in\dominant{\rankg}$
(\cite{AS}).

For complex numbers $a,b$
such that $b-a+1=\rankh$,
let $\C_{[a,b]}=\C {\bf 1}_{[a,b]}$ denote the one-dimensional
representation of $\Hecke{\rankh}$,
defined by
\begin{align}
\label{eq:H_one_dimensinal}
s_i\cdot {\bf 1}_{[a,b]}={\bf 1}_{[a,b]}
\quad \left(i=1,\dots,\rankh-1\right),\\
\e{i}\cdot {\bf 1}_{[a,b]}=
(a+i-1){\bf 1}_{[a,b]}\quad
\left(i=1,\dots,\rankh\right).
\end{align}

Let $\Wt( X)$
denote the space of weights of a $\lie{\rankg}$-module $X$.
Then,
clearly
 $$\Wt( (\C^{\rankg})^{\tensor \rankh}))
=\{
\sum_{i=1}^{\rankg}n_i \e{i}\mid
n_i\in \{
0,1,\dots,\rankh\},\,
\ \sum_{i=1}^{\rankg}n_i=\rankh\}.$$
For $\lambda\in \dual{\cartan{\rankg}}$,
let
\begin{align}
S(\lambda;\rankh):=
\left\{\mu
\in \dual{\cartan{\rankg}}\
\mid \lambda-\mu\in \Wt( (\C^{\rankg})^{\tensor \rankh}))
\right\}.
\end{align}
For $\lambda\in \dual{\cartan{\rankg}}$
and $\mu\in S(\lambda,\rankh)$,
define an $\Hecke{\rankh}$-module
\begin{equation}
\label{df:standard}
 {\standard}(\lambda,\mu):=
\Hecke{\rankh}\tensor_{({\mathcal H}_{\rankh_1}\tensor
\dots\tensor{\mathcal H}_{{\rankh}_{\rankg}})}
(\C_{[\mu_1,\lambda_1-1]}\tensor\dots\tensor
\C_{[\mu_{\rankg},\lambda_{\rankg}-1]}),
\end{equation}
where 
$\lambda_i=\lambda(\e{i})$,
$\mu_i=\mu(\e{i})$,
$\rankh_i=\lambda_i-\mu_i$.
Here
$\Hecke{{\rankh}_1}\otimes \Hecke{{\rankh}_2}
\tensor
\dots \tensor \Hecke{{\rankh}_{\rankg}}$
 is 
regarded as
a subalgebra of $\Hecke{{\rankh}}$
via the embeddings
$\Hecke{\rankh_k}\hookrightarrow \Hecke{\rankh}$
defined by
$\e{a} \mapsto  \e{a+\sum_{j=1}^{k-1}\rankh_j}$ and
$s_{a} \mapsto  s_{a+\sum_{j=1}^{k-1}\rankh_j}$.

Let $${\bf 1}_{\lambda,\mu}:=
{\bf 1}_{[\mu_1,\lambda_1-1]}\tensor
\dots \tensor {\bf 1}_{[\mu_{\rankg},\lambda_{\rankg}-1]}
\in {\standard}(\lambda,\mu).$$
Then the correspondence 
${\bf 1}_{\lambda,\mu}\mapsto 1$ defines an isomorphism
of
  $W_{\rankh}$-modules
\begin{equation}
\label{eq:standard_as_w}
\standard(\lambda,\mu)
\cong
 \C[W_{\rankh}/W_{\rankh_1}\times
\dots\times W_{\rankh_{\rankg}}].
\end{equation}

For a partition $\nu$ of $\rankh$,
let $U(\nu)$ denote the simple $W_\rankh$-module
associated with  $\nu$.
For $\lambda\in\dual{\cartan{\rankg}}$
and $\mu \in S(\lambda;\rankh)$,
let $\nu_{\lambda,\mu}$ 
denote the partition of $\rankh$
obtained
by forgetting the
order of $(\ell_1,\dots,\ell_{\rankh})$,
where $\ell_i=(\lambda-\mu)(\e{i})$.
Then by \eqref{eq:standard_as_w},
$\standard(\lambda,\mu)$ decomposes as 
\begin{align}
\label{eq:decomposition_of_standard}
\standard(\lambda,\mu)
\cong
U(\nu_{\lambda,\nu})
\oplus \bigoplus_{\nu >\nu_{\lambda,\nu}}U(\nu)^{\oplus c_{\nu}},
\end{align}
as a $W_{\rankh}$-module,
where $>$ is the dominance order in the set of the partitions
and $c_{\nu}$ is  some nonnegative integer
(see \cite{FH}, for example).
It is known that if $\lambda\in
\dominant{\rankg}$,
the $W_{\rankh}$-simple 
component $U(\nu_{\lambda,\mu})$ generates 
$\standard(\lambda,\mu)$ over $\Hecke{\rankh}$,
hence it has an unique simple quotient
$\simple(\lambda,\mu)$
which
contains 
$U(\nu_{\lambda,\mu})$ with multiplicity one
(\cite{Ze1,Ro}, see also \cite{S}).
The module $\standard(\lambda,\mu)$ with
$\lambda\in\dominant{\rankg}$
and $\mu\in S(\lambda;\rankh)$
is called a \emph{standard module} of $\Hecke{\rankh}$.

Let $W_{\lambda}\subset W_{\rankg}$
denote the stabilizer of $\lambda\in \dual{\cartan{\rankg}}$.
Notice that if $\mu\in S(\lambda;\rankh)$,
then $w\cdot \mu\in S(\lambda;\rankh)$\
for all $w\in W_{\lambda}$.
One has
\begin{equation}
\label{eq:iso_standard}
\begin{array}{rcl}
\standard(\lambda,\mu)
\cong
\standard(\lambda,\mu')
&\iff &\simple(\lambda,\mu)
\cong\simple(\lambda,\mu')\\
&\iff&
\mu'=w\cdot \mu
\text{ for some }
w\in W_{\lambda}
\end{array}
\end{equation}
for $\lambda\in\dominant{\rankg}$
and $\mu,\,\mu'\in S(\lambda;\rankh)$
(\cite{Ze1,Ro}).
It is known that any simple $\Hecke{\rankh}$-module is isomorphic to
$\simple(\lambda,w\cdot \mu)$
for some 
$\lambda,\,\mu\in \dominant{\rankg}$ and
$w\in W_{\lambda}\backslash W_{\rankg} /W_{\mu}$
such that 
$w\cdot \mu\in S(\lambda;\rankh)$
for some $\rankg\in \nnumber$.
(\cite{Ze1,Ro}).

Let $M(\lambda)$ be the Verma module of 
$\lie{\rankg}$ with highest weight $\lambda-\rho$ and
let $L(\lambda)$ denote its unique simple quotient.
Let $\lambda, \,\mu\in \dominant{\rankg}$ and
$w\in W_{\lambda}\backslash W_{\rankg} /W_{\mu}$ such 
that $w\cdot \mu\in S(\lambda;\rankh)$.
Then one of the main results in \cite{AS,S} is stated as follows:
\begin{align}
&F_{\lambda}(M(w\cdot \mu))\cong 
\standard(\lambda,w\cdot \mu),
\label{eq:image_of_verma_F}\\
%
&F_{\lambda}(L(w\cdot \mu))
\cong
\begin{cases}
\simple(\lambda,w\cdot\mu) &\text{if }
w\cdot \mu=w_{LR}\cdot \mu,\\
0&\text{otherwise,}
\end{cases}
\label{eq:image_of_simple_F}
\end{align}
where $w_{LR}$ denotes the longest length representative
of $w$ in the
double coset $W_{\lambda}w W_{\mu}$.

\section{Drinfeld Functor}

\subsection{Drinfeld Functor}
For a 
left $\Hecke{\rankh}$-module $\Hmodule$,
consider an $\Hecke{\rankh}\tensor 
\UN$-module  $\Hmodule\tensor (\Vbox)^{\tensor \rankh}$.
Here we regard $\C^{\ranky}$  as the vector
representation of $\lie{\ranky}$.
For $x\in \lie{\ranky}$
and $i=1,\dots,\rankh$,
let $\tau_i(x)$ denote the endomorphism
of $\Hmodule\tensor (\Vbox)^{\tensor \rankh}$
which acts as 
$x\in \lie{\ranky}$ on the 
$i$-th factor of the tensor product $(\Vbox)^{\tensor \rankh}$ 
and by identity on all the other factors.

Define an action of the Yangian $\YN$
on $\Hmodule\tensor (\Vbox)^{\tensor \rankh}$
by
\begin{align}
\label{Yaction}
\pi:t(u)\longmapsto &T_1(u-\e{1})T_2(u-\e{2})\dots
T_{\rankh}(u-\e{\rankh}),
\end{align}
where
$$T_i(u-\e{i})=1+\frac{1}{ u-\e{i}}\Id{i}$$
and $\Id{i}=\sum_{1\leq a,b\leq \ranky} 
\tau_i(\Mat{ab}{\ranky})\otimes \Mat{ab}{\ranky}
\in \End(\Hmodule\tensor (\Vbox)^{\tensor \rankh})
\tensor \End(\Vbox)$.
By the fact that
$S(\cartan{{\rankh}})$ is commutative, 
\eqref{Yaction} gives a
 well-defined action of $\YN$ 
(recall \eqref{eq:coproduct} and \eqref{evaluation}).

The symmetric group $W_{\rankh}$ naturally acts on 
$\Hmodule\tensor (\Vbox)^{\tensor {\rankh}}$
by $s_{ij}\mapsto K_{ij}P_{ij}$,
where $K_{ij}$ denotes its action on $\Hmodule$
and $P_{ij}$ denotes its action on
$(\Vbox)^{\tensor {\rankh}}$ by permutation.
Now define
\begin{align}
\label{Dfunctor1}
\Dfun{\rankh}(\Hmodule):=(\Hmodule\tensor (\Vbox)^{\tensor {\rankh}})
/\sum_{i=1}^{\rankh}{\rm Im}(s_i+1),
\end{align}
where ${\rm Im}(s_i+1)$ denotes the subspace
$(s_i +1)(\Hmodule\tensor (\Vbox)^{\tensor {\rankh}})$.
Let $[m\tensor u]$ denote the equivalence class
of $m\tensor u\in \Hmodule\tensor (\Vbox)^{\tensor {\rankh}}$
in $\Dfun{\rankh}(\Hmodule)$ .

The following proposition
is due to Drinfeld (\cite{D2}, see also
\cite{Ber,CP1}).
\begin{Pro}
\label{Pro:Drinfeld_fun}
Let $M$ be an  $\Hecke{\rankh}$-module.
Then,
$\pi$ induces an
 action  of
$\YN$  on the space
$\Dfun{\rankh}(M)$. 
\end{Pro}
\begin{proof}
It is enough to show that 
$\left(\prod_{i=1}^{\rankh}(u-\e{i})\right)\cdot  
\pi(t(u))
$
preserves the denominator space
of 
(\ref{Dfunctor1})
since $\left(\prod_{i=1}^{\rankh}(u-\e{i})\right)
\in\Center(\Hecke{\rankh})[u]$.
This  follows from the
  formula
$$\begin{array}{l}
(u-\e{i}+\Id{i})(u-\e{i+1}+\Id{i+1})s_i\\
\equiv s_i(u-\e{i}+\Id{i})(u-\e{i+1}+\Id{i+1})
+\left(s_i+1\right)[\Id{i+1},\Id{i}]
\end{array}
$$
($s_i=K_{i i+1}P_{ii+1}$),
which can be proven by direct calculations
using the defining relations \eqref{eq:h_def}
and the commutation relations
$[P_i,I_i]=[I_{i+1},I_{i}]$.
\end{proof}
The action of $\YN$ on the space $\Dfun{\rankh}(M)$
will be denoted by the same symbol $\pi$.

\medskip

Let $\YNcat$ and $\YSNcat$ denote
the category of finite-dimensional
representations of  $\YN$ and $\YSN$
respectively.
Then  $\Dfun{\rankh}$ defines an exact functor from 
$\Hcat$ to $\YNcat$ or $\YSNcat$.
The functor $\Dfun{\rankh }$ is called 
the {\emph
{Drinfeld Functor}} (\cite{D2}).
Note that our definition of the Drinfeld functor 
slightly differs from that of 
\cite{D2}.
A thorem of Drinfeld states that
if $\rankh<\ranky$,
$\Dfun{\rankh}$ gives a categorical equivalence
between $\Hcat$ and the certain subcategory
of $\YSNcat$ (\cite{D2}).
One can deduce its unpublished proof from
the paper of Chari-Pressly \cite{CP1},
in which the  categorical equivalnece
was generalized to the case
between the affine Hecke algebra and the quantum affine algebra.
However,
the  method in \cite{CP1} does not apply to the case
 $\rankh\geq \ranky$.

\subsection{}
The following proposition follows from
the  Frobenius-Schur duality.
\begin{Pro}
\label{y_decom_as_g}
Let $M$ be an
$ \Hecke{\rankh}$-module.
 Let 
$M= \bigoplus\limits_{\nu}U(\nu)^{\oplus c_{\nu}}$
$\left( c_{\nu}\in \nonnegative\right)$
be a decomposition of $M$ as a
$W_{\rankh}$-module.
Then,
$$\Dfun{\rankh}(M)
\cong\bigoplus\limits_\uesita{\nu}{\nu(\e{1})
\leq \ranky}L(\nu'+\rho)^{\oplus c_{\nu}}$$
 as a $\lie{\ranky}$-module,
where
$\nu'$ is the transpose of a partition $\nu$
identified with the dominant integral weight
of $\lie{\ranky}$.
\end{Pro}

See \cite{CP1}
for the proof of the
following proposition.
\begin{Pro} 
\label{pr:D_tensor}
Let $M_1$ and $M_2$ be representations of
$\Hecke{\rankh_1}$ and $\Hecke{\rankh_2}$
respectively.
Then,
$$\Dfun{\rankh_1}(M_1)\tensor \Dfun{\rankh_2}(M_2)
\cong \Dfun{\rankh_1+\rankh_2}
\left(\Hecke{\rankh_1+\rankh_2}\tensor_{
(\Hecke{\rankh_1}\tensor\Hecke{\rankh_2})
}
(M_1\tensor M_2)\right)$$
as a $\YN$ and $\YSN$-module.
\end{Pro}

The following formula 
was stated in \cite{Ber} as a conjecture.
\begin{Pro}
\label{pro:douti}
Let $M$ be an
$\Hecke{\rankh}$-module.
Then,
\begin{equation}
\label{eq:douti}
\pi(t_{ab}(u))\equiv
\delta_{ab}+\sum_{i=1}^{\rankh}\dfrac{1}{u-y_i}\tensor
\tau_i(\Mat{ab}{\ranky}).
\end{equation}
on the space $D_\rankh(M)$.
In particular,
$\pi(t_{ab}^{(d)})$ acts as
$\sum_{i=1}^{\rankh}y_i^d \tensor\tau_i(E_{ab})$.
\end{Pro}
\begin{proof}
We prove by induction on $k$ that
\begin{align*}
\left(1+\frac{\Id{1}}{u-\e{1}}
\right)\left(1+\frac{\Id{2}}{u-\e{2}}
\right)
\dots 
\left(1+\frac{\Id{k}}{u-\e{k}}
\right)
\equiv
\left(1+\sum_{i=1}^{k}\frac{1}{u-y_i}
\Id{i}\right)
\end{align*}
on $\Dfun{\rankh}(M)$.

There is nothing to prove  for $k=1$.
Let $k>1$.
By induction hypothesis,
\begin{align*}
&
\left(1+\frac{\Id{1}}{u-\e{1}}
\right)\left(1+\frac{\Id{2}}{u-\e{2}}
\right)
\dots 
\left(1+\frac{\Id{k}}{u-\e{k}}
\right)\\
&\equiv 1+
\left(\sum_{i=1}^{k-1}\frac{1}{u-y_i}\Id{i}\right)
+\frac{1}{u-\e{k}}\Id{k}+
\left(\sum_{i=1}^{k-1}
\frac{1}{u-y_i}\right)\cdot
\frac{1}{u-\e{k}}\Id{i} \Id{k} 
\end{align*}  
Since $\Id{i}\cdot \Id{k}=P_{ik}
\cdot \Id{k}$ and $P_{ik}\equiv -K_{ik}$
on $\Dfun{\rankh}(M)$,
\begin{align*}
&\frac{1}{u-\e{k}}\Id{k}+
\left(\sum_{i=1}^{k-1}
\frac{1}{u-y_i}\right)\cdot
\frac{1}{u-\e{k}}\Id{i} \Id{k} \\
&
\equiv
\left(
\frac{1}{u-\e{k}}
-\sum_{i=1}^{k-1}
K_{ik}\cdot
\frac{1}{u-y_i}\cdot
\frac{1}{u-\e{k}}
\right)\Id{k}\\
&=
\frac{1}{u-y_{k}}
\left(
u-y_{k}-
\sum_{i=1}^{k-1}
K_{ik}
\right)
\frac{1}{u-\e{k}}
\Id{k}
\\
&=
\frac{1}{u-y_{k}}
\cdot \left(u-\e{k}
\right)
\cdot 
\frac{1}{u-\e{k}}
\Id{k}=\frac{1}{u-y_{k}}\Id{k}.
\end{align*}
\end{proof}

Let
$\Lambda_{k}=\sum_{i=1}^{k}\e{i}\in
P_{\ranky}^+$
for $k=0,\dots,\ranky$
and let
$v_{\Lambda_k}$ denote the highest weight vector
of the simple $\lie{\ranky}$-module $L(\Lambda_k+\rho)$.
Then we can identify $v_{\Lambda_k}$ with
the  highest weight vector of 
the simple $\YN$-module ${\rm {ev}}_a
\left(L(\Lambda_{k}+\rho)\right)$
($a\in\C$).
It can be  checked directly that
 its
 weight
$(\zeta_1(u),\dots,\zeta_\ranky(u))$
is given by
\begin{align}
\label{eq:weight_of_ev}
\zeta_i(u)=
\begin{cases}
1+\frac{1}{u-i-a}&\text{if $1\leq i\leq k$}\\
1& \text{otherwise.}
\end{cases}
\end{align}

The following proposition,
which is easily follows from  Proposition \ref{pro:douti},
is
due to
Chari-Pressley \cite{CP1}.
\begin{Pro}
\label{Pro:image_of_one_dim}
Let $a,\,b$ 
    be complex numbers 
such that $b-a+1=\rankh$. Then,
as a $\YN$-module,
\begin{align}
\label{eq:iso_one-di}
\Dfun{\rankh}(\C_{[a,b]})
\cong
\begin{cases}
 {\rm {ev}}_a
\left(L(\Lambda_{\rankh}+\rho)\right)&\text{if }\rankh\leq \ranky,\\
0&\text{otherwise.}
\end{cases}
\end{align}
\end{Pro}
\medskip

\section{Main results}
\label{Main results}
\subsection{}
For a subspace $M'$ of an $\Hecke{\rankh}$-module $M$,
let $\Dfun{\rankh}(M')$ denote the image of 
$M'$ by the Drinfeld functor in $\Dfun{\rankh}(M)$.
The proof of the following proposition
is in Section
\ref{sec:proof_of_belong}.
\begin{Pro}
\label{pro:crucial}
Let $M$ be an $\Hecke{\rankh}$-module such that
$M$ is generated by some simple $W_{\rankh}$-submodule
$U$  of $M$.
Suppose that 
$\Dfun{\rankh}(U)$ is nonzero. 
Then,
$\Dfun{\rankh}(U)$ generates
$\Dfun{\rankh}(M)$ over $\YN$.
\end{Pro}

\medskip
Now 
let $\rankg\in\nnumber$.
For $\lambda\in \dual{\cartan{\rankg}}$,
let $${S}^{(\ranky)}(\lambda)=
\{\mu\in \cartan{\rankg}\mid
(\lambda-\mu)(\e{i})\in \{0,1,\dots,\ranky\}
\text{ for }i=1,\dots,\rankg\}.$$
For $\lambda\in \dual{\cartan{\rankg}}$
and $\mu\in S^{(\ranky)}(\lambda)$,
define a tensor product module $\Ystandard(\lambda,\mu)$ of $\YN$ by
\begin{equation}
\label{df:ystandard}
\Ystandard(\lambda,\mu):=\ev{\mu_1}(L(\Lambda_{\rankh_1}+\rho))
\otimes \dots \tensor\ev{\mu_{\rankg}}
(L(\Lambda_{\rankh_{\rankg}}+\rho)), 
\end{equation}
where $\mu_i=\mu(\e{i})$,
$\rankh_i=(\lambda-\mu)(\e{i})$.
Here $\YN$ acts via the coproduct
\eqref{eq:coproduct}.
Let $v_{\lambda,\mu}:=v_{\Lambda_{\rankh_1}}\tensor
\dots\tensor v_{\Lambda_{\rankh_{\rankg}}}\in
\Ystandard(\lambda,\mu)$.
Then by \eqref{eq:weight_of_ev},
$t_{ii}(u)\cdot v_{\lambda,\mu}=
\zeta_{\lambda,\mu;i}(u)v_{\lambda,\mu}$
for $i=1,\dots,\ranky$,
where
\begin{align}
\label{eq:highest_weight}
\zeta_{\lambda,\mu;i}(u)=\prod_\uesita{j=1,\dots,\rankh}{\rankh_j\geq i}
(1+\frac{1}{u-i-\mu_j}).
\end{align}

Let
$S^{(\ranky)}(\lambda;\rankh)=S(\lambda;\rankh)\cap S^{(\ranky)}(\lambda)$.
Notice that
for $\mu\in S(\lambda;\rankh)$,
the condition
$\mu\in S^{(\ranky)}(\lambda;\rankh)$
is equivalent to $\nu_{\lambda,\mu}(\e{1})\leq \ranky$ ,
where the partition $\nu_{\lambda,\mu}$
is identified with a dominant integral weight
(recall subsection \ref{subsection:Hecke_rep}).

\begin{Th}
\label{Th:properties_of_Dfun_standard} 
$ $

 $(1)$
The Drinfeld functor sends a standard module
of $\Hecke{\rankh}$ 
to zero or a highest weight module
of $\YN$.

$(2)$ More precisely,
let $\lambda
\in \dominant{\rankg}$
and 
 $\mu
\in {S}(\lambda;\rankh)$.
Then,
\begin{align}
\label{eq:explit_image_of_standard}
\Dfun{\rankh}(\standard(\lambda,\mu))
\cong
\begin{cases}
\Ystandard(\lambda,\mu)&\text{if } \mu\in S^{(\ranky)}(\lambda;\rankh),\\
0& \text{otherwise}.
\end{cases}
\end{align}
In particular, $\Ystandard(\lambda,\mu)$ is highest weight
with the highest weight vector $v_{\lambda,\mu}$
and the  highest weight
$\zeta_{\lambda,\mu}(u)$.
\end{Th}

\begin{proof}
(1)
Let $M$ be a standard module and
suppose that $\Dfun{\rankh}(M)\ne 0$.
 Since $M$ is a standard module,
$M\cong U(\nu)\oplus
 \bigoplus\limits_\uesita{\gamma>\nu}{|\gamma|=\rankh}
 U(\gamma)^{\oplus c_{\gamma}}$
and
$M=\Hecke{\rankh}\cdot U(\nu)$ 
for some partition $\nu$ of $\rankh$
such that $\nu(\e{1})\leq \ranky$.
By Proposition \ref{y_decom_as_g},
\begin{equation*}
\Dfun{\rankh}(M)
\cong L(\nu'+\rho)\oplus \bigoplus\limits_\uesita{\gamma'<\nu',\gamma(\e{1})
\leq \ranky}{|\gamma|=\rankh}L(\gamma'+\rho)^{\oplus c_{\gamma}}
\end{equation*}
But by  Proposition \ref{pro:crucial},
the highest weight vector of $L(\nu'+\rho)$ generates
$\Dfun{\rankh}(M)$
over $\YN$
in this decomposition.
Since the other $\lie{\ranky}$-weights appearing in 
$\Dfun{\rankh}(M)$ is smaller than $\nu'$
with respect to the dominance order,
it follows that $\Dfun{\rankh}(M)$
is a highest weight module
whose 
 highest weight vector is
the $\lie{\ranky}$-highest weight vector of
$L(\nu'+\rho)\subset \Dfun{\rankh}(M)$.
(2)
The isomorphism  
\eqref{eq:explit_image_of_standard} follows from
 Proposition \ref{pr:D_tensor}
and Proposition  \ref{Pro:image_of_one_dim}.
In fact,
\eqref{eq:explit_image_of_standard} holds without restriction
$\lambda\in\dominant{\rankh}$.
The rest of the statement follows from (1).
\end{proof}

\begin{Rem}
The fact that
$\Ystandard(\lambda,\mu)$ is highest weight
 for $\lambda\in\dominant{\rankg}$
was proved by Akasaka-Kashiwara (\cite{Ak})
in the case of the quantum affine algebra 
and by Nazarov-Tarasov (\cite{Tara1}) in the  case of the Yangian.
The above thorem provides  another  proof of it.
\end{Rem}
Let us call those $\YN$-modules
$M(\lambda,\mu)$ with
$\lambda\in \dominant{\rankg}$ and $\mu\in S^{(\ranky)}(\lambda;\rankh)$
{\it standard tensor product modules} of $\YN$.
By Theorem \ref{Th:properties_of_Dfun_standard}  (2),
 a standard tensor product  module
$\Ystandard(\lambda,\mu)$ has a 
 unique simple quotient, which is denoted by
$\Ysimple(\lambda,\mu)$. 
Then by \eqref{eq:highest_weight},
its Drinfeld polynomials
$Q_{\lambda,\mu}(u)=
(Q_{\lambda,\mu;1}(u),\dots,Q_{\lambda,\mu;\ranky-1}(u))$
are caluculated  as
\begin{equation}
\label{eq:doly_poly}
Q_{\lambda,\mu;k}(u)=\prod_\uesita{i=1,\dots,\rankg}{\lambda_i-\mu_i=k}
(u-\lambda_i),
\end{equation}
where $\lambda_i=\lambda(\e{i})$ and $\mu_i=\mu(\e{i})$
(\cite{Ak,CP1,Tara1}).

\begin{Th}
\label{Th:properties_of_Dfun_simple} 
$ $

$(1)$
The Drinfeld functor sends a simple
$\Hecke{\rankh}$-module to 
zero or a simple 
$\YN$-module.

$(2)$
More precisely,
let $\lambda
\in {\dominant{\rankg}}$ and
 $\mu
\in {S}(\lambda;\rankh)$.
Then,
\begin{align*}
\Dfun{\rankh}(\simple(\lambda,\mu))
\cong
\begin{cases}
\Ysimple(\lambda,\mu)&\text{if} \,\mu\in S^{(\ranky)}(\lambda;\rankh),\\
0& \text{otherwise}.
\end{cases}
\end{align*}

\end{Th}
\begin{proof}
(1)
Let $\simple$ be a simple
$\Hecke{\rankh}$-module.
Suppose that
$\Dfun{\rankh}(\simple)\ne 0$.
Let $V$ be a proper $\YN$-submodule 
of $\Dfun{\rankh}(\simple)$.
We suppose that $V\ne 0$ and deduce a contradiction.
Let $L$ be a simple $\lie{\ranky}$-submodule
of $V$.
Then $L=\Dfun{\rankh}(U)$ 
for some 
simple $W_{\rankh}$-submodule
$U$ of $M$.
But since $M$ is simple,
$M=\Hecke{\rankh}\cdot U$.
Thus by Proposition \ref{pro:crucial},
 $\Dfun{\rankh}(M)=\YN\cdot L\subset
\YN\cdot V$, which contradicts the assumption
that $V$ is proper.
(2)
follows from (1),
Theorem \ref{Th:properties_of_Dfun_standard} (2) and
and the fact that $\simple(\lambda,\mu)$
contains 
the simple $W_{\rankh}$-module $U(\nu_{\lambda,\mu})$ with multiplicity one.
\end{proof}

\begin{Rem}
By \eqref{eq:doly_poly},
it is easy to see that
every  simple $\YSN$-module
appears as the image of 
a simple $\Hecke{\rankh}$-module
by the Drinfeld functor for some $\rankh$. 
\end{Rem}

\begin{Rem}
When
$\lambda-\rho$ and $\mu-\rho$ 
are both dominant  weights,
$\Ysimple(\lambda,\mu)$ belongs to the class
of representations called tame
 (\cite{Tara2}).
Conversely,
any simple tame module is isomorphic to
$\Ysimple(\lambda,\mu)$ with $\lambda-\rho,\,\mu-\rho
\in\dominant{\rankg}$ for some $\rankh$.
\end{Rem}

\subsection{Multiplicity Formula}
For $\lambda,\mu\in P_{\rankg}^+$,
define
$$W_{\rankg}^{(\ranky)}(\lambda,\mu)=\{w\in
W_{\rankg}\mid
w\cdot \mu\in {S}^{(\ranky)}(\lambda;\rankh)\}\subset W_{\rankg}.$$
Then  $S^{(\ranky)}(\lambda;\rankh)=\bigsqcup_{\mu\in P_{\rankg}^+}
\bigcup_{w\in W_{\rankg}^{(\ranky)}(\lambda,\mu)}\{w\cdot \mu\}$.
Notice that
if $\mu\in S^{(\ranky)}(\lambda;\rankh)$,
then $w\cdot \mu\in S^{(\ranky)}(\lambda;\rankh)$
for all $w\in W_{\lambda}$.

\begin{Lem}
\label{Rem:same_poly}
Let $\lambda,\,\mu\in \dominant{\rankg}$
and $w,w'\in W_{\rankg}^{(\ranky)}(\lambda,\mu)$.
Then
$Q_{\lambda,w\cdot \mu}(u)=Q_{\lambda,w'\cdot\mu}(u)$
if and only if 
$w\equiv w'$ in the double coset
$ W_{\lambda}\backslash W_{\rankg}/ W_{\mu}$.
\end{Lem}
\begin{proof}
First notice that
the condition
$w\equiv w'$ in the double coset
$ W_{\lambda}\backslash W_{\rankg}/ W_{\mu}$
is equivalent to
the condition that
 the following sets
of pairs of complex numbers
are equal;
\begin{equation}
\label{temp_1}
\left\{(\lambda(\e{i}),w\cdot \mu(\e{i}))\mid i=1,\dots,\rankg
\right\},\quad
\left\{(\lambda(\e{i}),w'\cdot \mu(\e{i}))\mid i=1,\dots,\rankg
\right\}.
\end{equation}
Hence
the direction $\Leftarrow $ is easy to see.
$\Rightarrow$.
Let
\begin{equation}
\label{temp5.2.9_0}
\{(a_1,b_1),\cdots,(a_k,b_k)\},
\quad
\{(a_1,b_1'),\cdots, (a_k,b_k')\}
\end{equation}
be the result of removing
the common pairs from 
\eqref{temp_1}
so that
$a_i\geq a_{j}$ if $a_i-a_j\in\integer$
for $1\leq i<j\leq k$.
We suppose that $k\geq 1$ and deduce a contradiction.
By the assumption
$Q_{\lambda,w\cdot\mu}(u)=
Q_{\lambda,w'\cdot \mu}(u)$,
the differences $a_i-b_i$
and $a_i-b_i'$ are all $0$ or $\ranky$.

Now by the assumption, $b_1\ne b_1'$.
We can assume $a_1=b_1+\ranky=b_1'$.
Since
$\{b_i\}_{i=1}^k=\{b_i'\}_{i=1}^k
$,
there exists $p$ such that $b_p=b_1'$.
Then 
$a_p=b_p(=a_1)$ since $a_p=b_p$ or $b_p+\ranky$
and $a_1\geq a_p$.
Hence $(a'_1,b_1)=(a_p,b_p)$,
which contradicts the assumption that
there is no common
pair between the two sets
in \eqref{temp5.2.9_0}.
\end{proof}

\begin{Pro}
\label{pro:same_pol}
Let $\lambda,\,\mu\in \dominant{\rankg}$
and $w,w'\in W_{\rankg}^{(\ranky)}(\lambda,\mu)$.
Then the following conditions are equivalent:

$(1)$ $\Ystandard(\lambda,w\cdot \mu)\cong 
\Ystandard(\lambda,w'\cdot \mu)$.

$(2)$ $\Ysimple(\lambda,w\cdot \mu)\cong 
\Ysimple(\lambda,w'\cdot \mu)$.

$(3)$ $w\equiv w'$ in the double coset
$ W_{\lambda}\backslash W_{\rankg}/ W_{\mu}$.
\end{Pro}
\begin{proof}
(1) $\Rightarrow$ (2) follows from 
Theorem \ref{Th:properties_of_Dfun_standard} and
the fact that a  simple quotient
of a highest weight module is unique.
(2) $\Rightarrow$ (3) follows
 from Lemma \ref{Rem:same_poly}
and
(3) $\Rightarrow$ (1)
follows from \eqref{eq:iso_standard}.
\end{proof}

Let $\overline{W_{\rankg}^{(\ranky)}(\lambda,\mu)}$ denote the image
of $W_{\rankg}^{(\ranky)}(\lambda,\mu)$  in
the double coset
$W_{\lambda}\backslash W_{\rankg}/
W_{\mu}$.
By Proposition \ref{pro:same_pol},
each correspondence
\begin{align*}
w\mapsto \Ysimple(\lambda,w\cdot\mu),\quad
w\mapsto \Ystandard(\lambda,w\cdot\mu)
\end{align*}
defines  an injective map
from the set $\overline{W_{\rankg}^{(\ranky)}(\lambda,\mu)}$
to the set of equivalence classes of
finite-dimensional $\YSN$-modules.

Let $\leq$ denote the Bruhat ordering in $W_{\rankg}$
and let $P_{w,x}(q)$ denote the Kazhdan-Lusztig
polynomial associated with the Weyl group $W_{\rankg}$
(\cite{KL}).

In the following theorem
we state the multiplicity formula
of $\YSN$.
For simplicity,
we only consider the essential cases
when the  roots of Drinfeld polynomials
are integers.
\begin{Th}
\label{th:main}
Let $\lambda,\,\mu\in \indominant{\rankg}$
and $w\in \overline{W_{\rankg}^{(\ranky)}(\lambda,\mu)}$.

$(1)$
The family $\{\Ysimple(\lambda,x\cdot\mu)\mid
x\in \overline{W_{\rankg}^{(\ranky)}(\lambda,\mu)},
\,x_{LR}\geq w_{LR}\}$ is  exactly the set of 
all simple $\YSN$-modules which appear
as the composition factors  of
the standard tensor product module
  $ \Ystandard(\lambda, w\cdot \mu)$.
Moreover, their multiplicities
are expressed as
$$[\Ystandard(\lambda, w\cdot \mu),\Ysimple(\lambda, x\cdot \mu)]
=P_{w_{LR},x_{LR}}(1)$$
for $x\in \overline{W_{\rankg}^{(\ranky)}(\lambda,\mu)}$.

$(2)$ Conversely,
the simple 
$\YSN$-module $\Ysimple(\lambda, w\cdot \mu)$
is expressed as
$$
[\Ysimple(\lambda, w\cdot \mu)]=
\sum_\uesita{x\in W_{\rankg}^{(\ranky)}(\lambda,\mu)}
{{x}\geq w_{LR}}
(-1)^{\ell(w_{LR}w_0)-
\ell(xw_0)}
P_{xw_0,w_{LR}w_0}(1)[\Ystandard(\lambda, x\cdot \mu)].
$$
in the Grothendieck group
of $\YSNcat$,
where $w_0$ denotes the longest
element of $W_{\rankg}$.
\end{Th}
\begin{proof}
Due to the well-known
Kazhdan-Lusztig conjecture (\cite{BB,BK}) and the 
translation principle (\cite{Ja}),
one has
\begin{align}
&[M(w\cdot\mu)]=
\sum_\uesita{x\in W/W_{\mu}}{x_R\geq w_R}
P_{w_R,x_R}(1)[L(x\cdot \mu)]
\label{multi_g_verma}\\
&[L(w\cdot \mu)]=
\sum_\uesita{x\in W_{\rankg}}{x\geq w_R}
(-1)^{\ell(w_R w_0)-\ell(xw_0)}
P_{xw_0,w_R w_0}(1)[M(x\cdot \mu)],
\label{multi_g_simple}
\end{align}
for $\mu\in P_{\rankg}^+$ and $w\in W/W_{\mu}$
  in the Grothendieck group of 
the category
$\BGG$  of $\lie{\rankg}$,
where $w_R$ denote the longest length representative in the
coset $w W_{\mu}$.

Then by
\eqref{eq:image_of_verma_F},
\eqref{eq:image_of_simple_F},
Theorem \ref{Th:properties_of_Dfun_standard}
and Theorem \ref{Th:properties_of_Dfun_simple},
applying the exact functor $\Dfun{\rankh}\circ
F_{\lambda}$,
one has
\begin{align}
&[\Ystandard(\lambda,w\cdot \mu)]=
\sum_\uesita{x\in \overline{
W_{\rankg}^{(\ranky)}(\lambda, \mu)}}{x_{LR}\geq w_R}
P_{w_R,x_{LR}}(1)[\Ysimple(\lambda,x\cdot \mu)]
\label{eq:multi_y_standard}\\
&[\Ysimple(\lambda,w\cdot \mu)]=
\sum_\uesita{x\in W_{\rankg}^{(\ranky)}(\lambda, \mu)}{x\geq w_{LR}}
(-1)^{\ell(w_{LR}w_0)-\ell(xw_0)}
P_{xw_0,w_{LR}w_0}(1)[\Ystandard(\lambda,x\cdot \mu)]
\label{eq:multi_y_simple}
\end{align}
in the Grothendieck
group of  $\YSNcat$.
Hence (2) is proved
and (1) follows form
Proposition \ref{pro:same_pol}
and
\eqref{eq:multi_y_standard}.
\end{proof}

\begin{Rem}
If $x_R\geq w$ (resp. $x_L\geq w$),
then 
$x_R\geq w_R$ and $P_{w,x_R}(q)=P_{w_R,x_R}(q)$
(resp. $ x_L\geq w_L$ and $P_{w,x_L}(q)=P_{w_L,x_L}(q)$),
where  $w_L$ denotes the longest length representative
in the coset $W_{\lambda}w$
(see \cite[Corollary 7.14]{H}, for example).
Hence it follows that $P_{w,x_{LR}}(q)=P_{w_R,x_{LR}}(q)=
P_{w_{LR},x_{LR}}(q)$.
\end{Rem}

\begin{Rem}
Let $\lambda-\rho$ and $\mu-\rho$ 
are both dominant integral weights,
one can obtain the
resolution 
of 
$V(\lambda,\mu)$ 
by applying the exact functor $\Dfun{\rankh}\circ
F_{\lambda}$.
This provides an alternative proof
of the resolutions
of the simple {\emph{elementary module}} (\cite{Tara1})
constructed by Cherednik (\cite{Ch1}).
A generalization of such resolutions  is possible
to some extent
by using the generalized BGG resolution
obtained in \cite{GJ}
(see \cite{S} for the corresponding statement
in the case of the degenerate affine Hecke algbera).
\end{Rem}

\section{Proof of Proposition \ref{pro:crucial}}
\label{sec:proof_of_belong}
Let $M$ be an $\Hecke{\rankh}$-module
and $U\subset M$ be a simple $W_{\rankh}$-module
such that $\Dfun{\rankh}(U)\ne 0$
as in the proposition.
The proof is divided into 4 parts.
\subsection{}
Let $\tilde \Hcat$ be the
category
of $\Hecke{\rankh}$-modules
who decompose into
(possibly infinite) direct sum
of finite-dimensional $W_{\rankh}$-modules.
Let $\tilde \YNcat$ be 
the
category
of $\YN$-modules
who decompose into
(possibly infinite) direct sum
of finite-dimensional $\lie{\ranky}$-modules.
Extend the Drinfeld functor $\Dfun{\rankh}$
to the functor 
from the category $\tilde \Hcat$
to the category $\tilde \YNcat$.
It is easy to see that the extended functor
$\Dfun{\rankh}$ is still 
an exact functor.

Let $\nu=\sum_{i=1}^{\ranky}\nu_i\e{i}$ be the partition such that
$U\cong U(\nu')$,
where $\nu'$ denotes the transpose 
of $\nu$ as before.
Define an $W_{\rankh}$-module
$$J_{W}(\nu):=
\C[W_{\rankh}]\tensor_{(\C[W_{\nu_1}]
\tensor\dots\tensor \C[W_{\nu_{\ranky}}])}( \C 
{\bf{1}}_{\sign,\nu_1}
\tensor\dots\tensor
\C{\bf{1}}_{\sign,\nu_{\ranky}}),$$
where $\C {\bf{1}}_{\sign,\nu_i}$ is 
one-dimensional representation of $\C[W_{\nu_i}]
$ such that $s_i\cdot{\bf{1}}_{\sign,\nu_i}=
-{\bf{1}}_{\sign,\nu_i}$.
It is well-known 
that there exists a surjective
homomorphism $\varphi_{W}:J_{W}(\nu)
\twoheadrightarrow 
U(\nu')$
of $W_{\rankh}$-modules.
Let $J(\nu)$ be  an $\Hecke{\rankh}$-module 
defined by
$$J(\nu):=\Hecke{\rankh}\tensor_{\C[W_{\rankh}]}
J_{W}(\nu).$$
Then $J(\nu)$ is an object of
$\tilde \Hcat$.
Let $\varphi:J(\nu)\twoheadrightarrow M$ be the 
surjective
$\Hecke{\rankh}$-homomorphism
induced by $\varphi_W$.
Then, the $\YN$-homomorphism
$$\Dfun{\rankh}(\varphi):\Dfun{\rankh}(J(\nu))\longrightarrow
\Dfun{\rankh}(M)$$
is surjective and $\Dfun{\rankh}(\varphi)(J_W(\nu))=
\Dfun{\rankh}(U)$.
Hence it is suffice to prove that
\begin{align}
\label{eq:problem}
\Dfun{\rankh}(J(\nu))=\YN\cdot \Dfun{\rankh}(J_W(\nu)).
\end{align}

\subsection{}
Recall the filterations
on $\Hecke{\rankh}$
and $\YN$.
For a $\gr \Hecke{\rankh}$-module $\bar M$,
let $\bar \Dfun{\rankh}(\bar M)$
denote the $\gr \YN$-module
$\bar \Dfun{\rankh}(\bar M):=
(\bar M\tensor (\Vbox)\tensorN)/\sum {\rm Im}(s_i+1)$,
in which $U(\lie{\ranky}[t])\cong \gr \YN$
acts as
$$E_{ij}\tensor t^r\mapsto
\sum_{a=1}^{\rankh}\bar y_a^r\tensor \tau_{a}(E_{ij}).$$
Then
$\bar \Dfun{\rankh}$ defines an exact functor
from the category
of $\gr \Hecke{\rankh}$-modules who decomposes
into
 direct sum
of finite-dimensional $W_{\rankh}$-modules
to the category of 
$U(\lie{\ranky}[t])$-modules
decomposes into
 direct sum
of finite-dimensional $\lie{\ranky}$-modules.

Now introduce a filteration on
$J(\nu)$
by the followings;
\begin{align*}
&F_{-1} J(\nu)=0,\,F_{0} J(\nu)=J_W(\nu),
\,F_i J(\nu)=\left(F_i \Hecke{\rankh}\right)\cdot
\left( F_0 J(\nu)\right)
\quad \left( i\in \nnumber\right).
\end{align*}
Let $\bar J(\nu)$ 
denote the 
the corresponding graded module of 
$J(\nu)$.
Then,
as a $\gr \Hecke{\rankh}$-module,
\begin{align*}
\bar J(\nu)&=\gr \Hecke{\rankh}\cdot J_W(\nu)\\
&\cong
(\C[\bar y_1,\dots,\bar y_{\rankh}]
\tensor \C[W_{\rankh}])\tensor_{\C[W_{\rankh}]}
J_W(\nu).
\end{align*}
In particular,
$\bar J(\nu)
\cong \C[\bar y_1,\dots,\bar y_{\rankh}]
\tensor_{\C} J_W(\nu)$
as a $\C$-vector  space.

Introduce a filteration on $\Dfun{\rankh}(J(\nu))$
induced from that of $J(\nu)$;
\begin{align}
\label{filter:D(J)}
F_i \Dfun{\rankh}(J(\nu)):=
\Dfun{\rankh}(F_i J(\nu)).
\end{align}
By Proposition \ref{pro:douti},
one can easily check that
$F_i \YN\cdot F_j \Dfun{\rankh}( J(\nu))
\subset F_{i+j}
\Dfun{\rankh}(J(\nu))$.
Let $\gr \Dfun{\rankh}(\bar J(\nu))$
denote the corresponding graded module of
 $\Dfun{\rankh}(J(\nu))$. 
Then,
by Proposition \ref{pro:douti} and
the fact that 
$\bar \Dfun{\rankh}$ is exact, 
one has
$$\gr \Dfun{\rankh}(J(\nu))
\cong \bar \Dfun{\rankh}(\bar J(\nu))$$
as a  $\gr\YN\cong U(\lie{\ranky}[t])$-module.
Now \eqref{eq:problem} is reduced
to  the following proposition.
\begin{Pro}
As a $U(\lie{\ranky}[t])$-module,
\label{Pro:problem2}
\begin{align*}
\bar \Dfun{\rankh}(\bar J(\nu))
=U(\lie{\ranky}[t])\cdot
\bar \Dfun{\rankh}(J_W(\nu)),
\end{align*}
where $\bar \Dfun{\rankh}(J_W(\nu))$ denotes
the image of $J_W(\nu)\subset \bar J(\nu)$
in $\bar \Dfun{\rankh}(\bar J(\nu))$.
\end{Pro}

\subsection{}
As a preparation for the proof of Proposition
\ref{Pro:problem2},
we shall first consider the simplest case
when
$\nu=\rankh \e{i}$ for some $i$.
Put $\bar J(\rankh)=\bar J(\rankh \e{i})$.
Then as a $\gr \Hecke{\rankh}$-module,
\begin{align*}
\bar J(\rankh)
\cong\C[\bar y_1,\dots,\bar y_{\rankh}]{\bf 1}_{\sign},
\end{align*}
where 
$W_{\rankh}$ act on the right-hand-side by
$w\cdot (f{\bf 1}_{\sign})=(-1)^{\ell(w)}w(f){\bf 1}_{\sign}$.

Let ${{\bf i}}=(i_1,\dots,i_{\ranky})$
be a permutation of $(1,2,\dots,\ranky)$. 
For $\gamma=\sum_{i=1}^{\ranky}\gamma_i \e{i}\in \Wt
\left((\Vbox)^{\tensor \rankh}\right)$,
let
$$u_{{{\bf i}}}(\gamma):=
u_{i_1}^{\tensor \gamma_{i_1}}\tensor
u_{i_2}^{\tensor \gamma_{i_2}}\tensor
\dots
\tensor
u_{i_{\ranky}}^{\tensor \gamma_{i_{\ranky}}}\in
(\Vbox)\tensorN.$$
\begin{Lem}
For a fixed permutation
${{\bf i}}=(i_1,\dots,i_{\ranky})$
of $(1,2,\dots,\ranky)$,
the set
\begin{align}
\label{eq:basis}
\left\{\left[\bar y_1^{d_1}\bar y_2^{d_2}
\dots \bar y_{\rankh}^{d_{\rankh}}{\bf 1}_{\sign}
\tensor u_{{{\bf i}}}(\gamma)\right]
\mid
\begin{array}{l}
\gamma=\sum_{i=1}^{\ranky}\gamma_i \e{i}\in \Wt
\left((\Vbox)^{\tensor \rankh}\right),\\
d_m\geq d_{m+1}
\text{ {\rm if} } \sum_{j=1}^{a-1}\gamma_{i_j}<m<\sum_{j=1}^{a}\gamma_{i_j}\\
\text{{\rm for some} }a.
\end{array}
\right\}
\end{align}
forms a $\C$-basis of 
$\bar \Dfun{\rankh}\left(\bar J(\rankh)\right)$.
\end{Lem}
\begin{proof}
Clearly
\begin{align*}
\bar \Dfun{\rankh}\left(\bar J(\rankh)\right)
\cong \bigoplus_{\gamma\in
\Wt((\Vbox)\tensorN)}
(\bar J(\rankh)\tensor \left[ (\Vbox)\tensorN\right]_{\gamma})
/\sum {\rm Im}(s_i+1)
\end{align*}
as a $\C$-vector space.
In the case when $\gamma=\rankh \e{i}$ for some $i$,
it is well-known that
\begin{align}
\label{eq:basis_J}
(\bar J(\rankh)\tensor \left[ (\Vbox)\tensorN \right]_{\rankh \e{i}})
/\sum {\rm Im}(s_i+1)
=\bigoplus_{d_1\geq \dots \geq d_{\rankh}} \C[
\bar y_1^{d_1}
\dots \bar y_{\rankh}^{d_{\rankh}}{\bf 1}_{\sign}
\tensor u_i^{\tensor \rankh}].
\end{align}
For a general weight $\gamma\in
\Wt((\Vbox)\tensorN)$,\
notice that
the correspondence
$w\mapsto w\cdot u_{{{\bf i}}}(\gamma)$
defines an isomorphism
$\C[W_{\rankh}/
W_{\gamma_{i_1}}\times
\dots \times W_{\gamma_{i_{\ranky}}}]
\cong \left[ (\Vbox)\tensorN\right]_{\gamma}
$.
On the other hand,
$\bar J(\rankh)\cong
\bar J(\gamma_{i_1})\tensor 
\dots\tensor \bar J(\gamma_{i_{\ranky}})$
as a $W_{\gamma_{i_1}}\times
\dots \times W_{\gamma_{i_{\ranky}}}$-module.
Hence, by using the same argument as in the Proposition \ref{pr:D_tensor},
the statement reduces to \eqref{eq:basis_J}.
\end{proof}

For
$k\in \{1,\dots,\ranky\}$,
let $U(\lie{\ranky}[t])_k$ denote the subalgebra of 
$U(\lie{\ranky}[t])$
generated by 
the elements $E_{ik}\tensor t^r$
$\left( i=1,\dots,\ranky,\,
r\in \nonnegative\right)$.
\begin{Lem}
\label{Pro:essential}
For any  $k\in \{1,\dots,\ranky\}$,
\begin{align}
\label{eq:essential}
\bar \Dfun{\rankh}\left(\bar J(\rankh)\right)=
U(\lie{\ranky}[t])_k\cdot \left[
{\bf 1}_{\sign}\tensor u_k^{\tensor \rankh}\right].
\end{align}
\end{Lem}
\begin{proof}
Let ${{\bf i}}=(i_1,\dots,i_{\ranky})
=(1,2,\dots,k-1,k+1,k+2.\dots,\ranky,k)$.
For $m=0,1,\dots,\rankh$,
let
$\bar \Dfun{\rankh}\left(\bar J(\rankh)\right)_{k,m}$
be the subspace of
$\bar \Dfun{\rankh}\left(\bar J(\rankh)\right)$
spanned by the
vectors of the form \eqref{eq:basis}
such that $d_{j}=0$ for all $j>m$.
We prove by inducition on $m$
that $\bar \Dfun{\rankh}\left(\bar J(\rankh)\right)_{k,m}
\subset U(\lie{\ranky}[t])_k\cdot \left[
{\bf 1}_{\sign}\tensor u_k^{\tensor \rankh}\right]$.

Let $m=0$.
Then for any $\gamma=\sum_{i=1}^{\ranky}\gamma_i\e{i}\in
\Wt((\Vbox)^{\tensor \rankh})$,
\begin{align*}
&(E_{i_1 k})^{\gamma_{i_1}}
(E_{i_2 k})^{\gamma_{i_2}}
\dots 
(E_{i_{\ranky-1} k})^{\gamma_{i_{\ranky-1}}}
\cdot \left[
{\bf 1}_{\sign}\tensor u_k^{\tensor \rankh}\right]\\
&\equiv \frac{\rankh !}{\gamma_k !}
[{\bf 1}_{\sign}\tensor
u_{{{\bf i}}}(\gamma)],
\end{align*}
where $E_{ij}=E_{ij}\tensor 1\in \lie{\ranky}[t]$.
Hence 
$\bar \Dfun{\rankh}\left(\bar J(\rankh)\right)_{k,0}
\subset U(\lie{\ranky}[t])_k\cdot \left[
{\bf 1}_{\sign}\tensor u_k^{\tensor \rankh}\right]$.

Next let $m>0$
and suppose that
$\bar \Dfun{\rankh}\left(\bar J(\rankh)\right)_{k,m-1}
\subset U(\lie{\ranky}[t])_k\cdot \left[
{\bf 1}_{\sign}\tensor u_k^{\tensor \rankh}\right]$.
Let $a$ be the integer such that
$\sum_{j=1}^{a-1}\gamma_{i_j}<m\leq \sum_{j=1}^{a}
\gamma_{i_j}$.
Then,
one has
\begin{align*}
&(E_{i_a k}\tensor t^{d_m})\cdot
[(\bar y_1^{d_1} \bar y_2^{d_2}\dots
\bar y_{{m-1}}^{d_{m-1}}{\bf 1}_{\sign})
\tensor
u_{{{\bf i}}}(\gamma-\e{i_a}+\e{k})]\\
&\equiv
(\rankh-m+1)[(\bar y_1^{d_1} \bar y_2^{d_2}\dots
\bar y_{{m}}^{d_{m}}{\bf 1}_{\sign})
\tensor
u_{{{\bf i}}}(\gamma)]
\end{align*}
for  $d_m>0$.
Here the equality holds modulo
 $\bar \Dfun{\rankh}\left(\bar J(\rankh)\right)_{k,m-1}$
if  $a=n$ (i.e, $i_a=k$).
Hence by inducition
hypothesis,
$\bar \Dfun{\rankh}\left(\bar J(\rankh)\right)_{k,m}
\subset U(\lie{\ranky}[t])_k\cdot \left[
{\bf 1}_{\sign}\tensor u_k^{\tensor \rankh}\right]$.
\end{proof}

\subsection{}
Let turn ourselves back to
the proof of Proposition \ref{Pro:problem2}.
The following lemma
is an analogue of Proposition \ref{pr:D_tensor}.
\begin{Lem}
As a $\lie{\ranky}[t]$-module,
$$\bar \Dfun{\rankh}\left(\bar J(\nu)\right)
\cong
\bar \Dfun{\nu_1}\left(\bar J(\nu_1)\right)
\tensor
\bar \Dfun{\nu_2}\left(\bar J(\nu_2)\right)
\tensor
\dots
\tensor
\bar \Dfun{\nu_{\ranky}}\left(\bar J(\nu_{\ranky})\right).
$$
\end{Lem}

Now let us complete the proof of Proposition
\ref{Pro:problem2}.
Let $v_i=\left[
{\bf 1}_{\sign}\tensor u_i^{\tensor \nu_i}\right]$
and $v_{[k]}:=v_{k}\tensor\dots\tensor v_{\ranky}$.
Notice that
$U(\lie{\ranky}[t])_k\cdot
v_{a}\subset
\C v_{a}$
if $a\ne k$.
Hence one can show by inducition on $k$
that
$$
U(\lie{\ranky}[t])_1\cdot
v_1\
\tensor
U(\lie{\ranky}[t])_2\cdot
v_2\tensor\dots
\tensor
U(\lie{\ranky}[t])_k\cdot
v_k\tensor
\C v_{[k+1]}
\subset
U(\lie{\ranky}[t])\cdot
v_{[1]}.
$$
Now Proposition
\ref{Pro:problem2}
follows from
 Lemma
\ref{Pro:essential}.
\begin{flushright}
{\it End of Proof}.
\end{flushright}

\bigskip

{\bf Acknowledgement }
I am deeply gratefull to Akihiro Tsuchiya
for his constant support,
including valuable advice
and carefully reading of the manuscript.
I am also gratefull to Tomoki Nakanishi
for his encouragement
and valuable advice.
I whould like to thank Takeshi Suzuki
for the collaboration in \cite{AS},
which lead to this work.

\bibliography{myref}
\end{document}